\documentclass[a4paper,11pt]{article}

\usepackage[utf8]{inputenc}
\usepackage[T1]{fontenc}
\usepackage{bbm}
\usepackage{amssymb,amsmath,amsthm}

\usepackage{mathtools}

\usepackage{authblk}

\usepackage{booktabs}

\usepackage{multirow}

\usepackage{biblatex}

\usepackage{tabularx}

\usepackage{xspace}

\usepackage{todonotes}

\usepackage{tikz}

\usepackage{tikzscale}

\usepackage{pgfplots}
\pgfplotsset{compat=1.14}
\usepgfplotslibrary{fillbetween}

\usepackage{algorithm2e}

\usepackage{siunitx}

\graphicspath{{Figures/}}

\usetikzlibrary{calc,intersections}

\tikzset{
  >=stealth,
  graphnode/.style={draw,circle, inner sep=1.},
  graphedge/.style={very thin},
  match/.style={dashed, thin}
}

\usepackage{subcaption}

\usepackage{hyperref}

\hypersetup {
  hidelinks,
}
\usepackage{cleveref}

\bibliography{references}

\newcommand{\R}{\mathbb{R}}
\newcommand{\N}{\mathbb{N}}
\newcommand{\Q}{\mathbb{Q}}

\newcommand{\Rp}{\R_{\geq 0}}

\theoremstyle{plain}
\newtheorem{thm}{Theorem}[section]
\newtheorem{lem}[thm]{Lemma}

\theoremstyle{definition}

\newcommand{\define}{\coloneqq}

\newcommand{\st}{\textrm{s.t.\xspace}}

\newcommand{\ie}{i.e.,\xspace}

\newcommand{\Series}{\mathcal{S}}

\renewcommand{\O}{\mathcal{O}}

\newcommand{\simple}{\mathrm{sim}}
\newcommand{\improved}{\mathrm{imp}}

\newcommand{\email}[1]{\href{mailto:#1}{\nolinkurl{#1}}}

\newcommand{\DTWMean}{\textsc{DTW-Mean}\xspace}

\renewcommand{\P}{\mathcal{P}}

\DeclareMathOperator{\lb}{lb}
\DeclareMathOperator{\ub}{ub}

\DeclareMathOperator{\NLP}{NLP}
\DeclareMathOperator{\proj}{proj}

\DeclareMathOperator{\dtw}{dtw}

\DeclareMathOperator{\conv}{conv}

\DeclarePairedDelimiter\floor{\lfloor}{\rfloor}

\theoremstyle{remark}

\newtheorem{remark}{Remark}

\DontPrintSemicolon
\RestyleAlgo{ruled}
\SetArgSty{}
\SetCommentSty{}

\SetKwComment{KwComment}{$\triangleright$}{}

\SetKwInOut{Input}{Input}
\SetKwInOut{Output}{Output}

\SetKw{Continue}{continue}
\SetKw{Break}{break}
\SetKw{KwBy}{by}
\SetKw{KwDownto}{downto}

\newbibmacro*{bbx:parunit}{%
  \ifbibliography
  {\setunit{\bibpagerefpunct}\newblock
    \usebibmacro{pageref}%
    \clearlist{pageref}%
    \setunit{\adddot\par\nobreak}}
  {}}

\renewbibmacro*{doi+eprint+url}{%
  \usebibmacro{bbx:parunit}
  \iftoggle{bbx:doi}
  {\printfield{doi}}
  {}%
  \iftoggle{bbx:eprint}
  {\usebibmacro{eprint}}
  {}%
  \iftoggle{bbx:url}
  {\usebibmacro{url+urldate}}
  {}}

\newcommand{\problemdef}[3]{
  \begin{center}
     \begin{minipage}{0.95\textwidth}
       \noindent
       \textsc{#1}
       \vspace{.1cm}

       \setlength{\tabcolsep}{3pt}
       \begin{tabularx}{\textwidth}{@{}lX@{}}
         \textbf{Input:} 		& #2 \\
         \textbf{Question:} 	& #3
       \end{tabularx}
     \end{minipage}
  \end{center}
}

\begin{document}

\title{Mathematical Programming Models for Mean Computation in Dynamic Time Warping Spaces}
\author[1]{Vincent Froese}
\author[2]{Christoph Hansknecht}

\affil[1]{\small
  Technische Universit\"at Berlin, Faculty~IV, Institute of Software Engineering and Theoretical Computer Science, Algorithmics and Computational Complexity,\protect\\
  \email{vincent.froese@tu-berlin.de}}
\affil[2]{\small Institute for Mathematical Optimization, TU Braunschweig,\protect\\ \email{c.hansknecht@tu-braunschweig.de}}

\date{\today}

\maketitle

\begin{abstract}
  The dynamic time warping (dtw) distance is an established tool for
  mining time series data. The \DTWMean problem consists of comptuing
  a series which minimizes the so-called Fréchet function, that is, the sum of
  squared dtw-distances to a given sample of time series. \DTWMean
  is NP-hard and intractable in practice. So far, this challenging
  problem has been solved by various heuristic approaches without any
  performance guarantees.

  We give a polynomial-time algorithm yielding lower bounds on
  the domain of a mean time series which translate into lower
  bounds on the Fréchet function.  We then formulate the problem as a
  discrete nonlinear optimization problem based on network flows. We introduce several mixed-integer nonlinear programming
  (MINLP) formulations in order to solve \DTWMean optimally.
  Our formulations are based on techniques such as
  outer approximations and nonlinear reformulations of the well-known
  big $M$ indicator constraints.

  Finally, we conduct several computational experiments to compare
  the different formulations on several instances derived from
  the UCR Time Series Classification Archive. While in general \DTWMean
  remains quite challenging, our fomrulations yield good
  results in several important specialized probem settings.

  \medskip


\end{abstract}

\section{Introduction}
\label{sec:introduction}

Dynamic time warping (dtw) is a widely used distance measure for distance-based time series mining~\cite{BLBLK17,AML19}.
It allows to cope with temporal variation in the data via nonlinear alignments between two input time series (see \Cref{sec:prelims} for details).
Averaging a sample of time series under the dtw-distance is a challenging optimization problem
in dtw-based time series mining.
Given samples~$s_1,\ldots,s_k$, the formal problem is to find a time series~$z$ with minimum \emph{Fr\'echet} variance
\[F(z) = \frac{1}{k}\sum_{i=1}^k\dtw(z,s_i)^2,\]
where~$\dtw(z,s_i)$ denotes the dtw-distance between~$z$ and~$s_i$.
A \emph{mean} is any time series minimizing~$F$.

\paragraph*{Related Work.}
It is known that a mean (of length at most~$nk$, where~$n$ is the maximum length of any input series) always exists~\cite{JS16}. Even for binary input series, the problem is known to be NP-hard, W[1]-hard with respect to the number~$k$ of samples (that is, presumably not solvable in~$f(k)\cdot n^{\O(1)}$ time for any function~$f$) and even not solvable in ~$n^{o(k)}\cdot f(k)$ time assuming the \emph{Exponential Time Hypothesis}\footnote{An assumption in complexity theory asserting that \textsc{3-SAT} is not solvable in time~$\O(2^{cn})$ for some constant~$c>0$, where~$n$ is the number of variables~\cite{IP01}.}~\cite{BFN18}.
The currently fastest exact algorithm uses dynamic programming and runs in~$\O(n^{2k+1}2^kk)$ time~\cite{BFFJNS19}.
Over the past decade, various heuristic approaches have been developed~\cite{HNF08,PKG11,CB17,SJ18,LZZ19}.
However, they all come without any theoretical performance guarantees and have been shown to yield
poor results in practice~\cite{BFFJNS19}.

\section{Preliminaries}
\label{sec:prelims}

\paragraph*{Time Series.}
For $n\in\N$, let~$[n]$ denote the set~$\{1,\ldots,n\}$.
We consider \emph{finite univariate rational time series}, which we will simply denote as
time series.  A time series is a sequence $s \in \Q^{n}$ for some
$n \in \N$.  We let $\Series_n$ be the set of all time series of length
$n$ and $\Series \define \bigcup_{n \in \N} \Series_n$ be the set of all time
series.

\paragraph*{Dynamic Time Warping.}
The \emph{diagonal grid graph} $D(m, n)$ for $m, n \in \N$ is the directed graph with
vertices
\[
  V(D(m,n)) \define [m]\times[n]
\]
and arcs
\[
  \begin{aligned}
    A(D(m,n)) \define & \{ ((i,j),(i+1, j+1)) \mid i \in [m-1], j \in [n-1] \} \: \cup \\
    & \{ ((i,j),(i+1, j)) \mid i \in [m-1], j \in [n] \} \: \cup \\
    & \{ ((i,j),(i,j+1)) \mid i \in [m], j \in [n-1] \}. \\
  \end{aligned}
\]
The \emph{origin} $s_{m, n}$ of $D(m,n)$ is defined as the vertex $(1, 1)$.
The set of \emph{destinations} is given as
\[
  T(D(m, n)) \define \{ (m, j) \mid j \in [n] \} \subseteq V(D(m, n)).
\]

The \emph{destination} of $D(m,n)$ is the vertex $t_{m, n} \define (m, n)$.
A \emph{warping path} of order $m \times n$ is an $s_{m, n}$-$t_{m,n}$-path through
$D(m,n)$. We let $\P_{m ,n}$ be the set of all warping paths of order~$m\times n$.

\begin{figure}
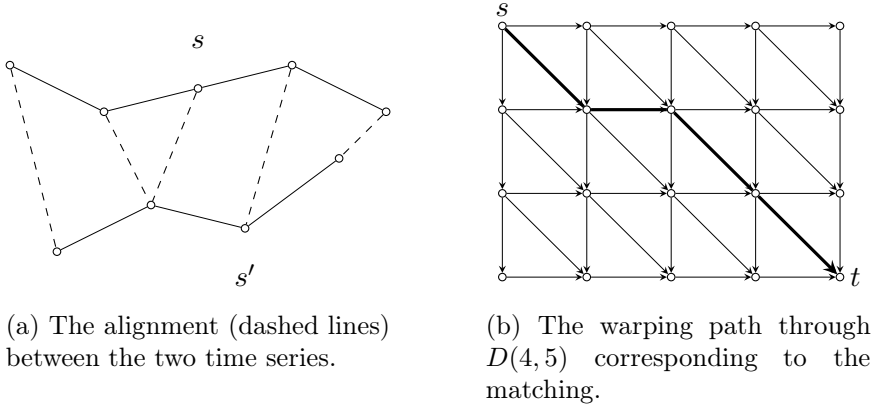

  \centering
  \begin{subfigure}[t]{.4\textwidth}
    \includegraphics[width=\textwidth]{Alignment}
    \caption{The alignment (dashed lines) between the two time series.}
  \end{subfigure}
  \hspace{1cm}
  \begin{subfigure}[t]{.4\textwidth}
    \includegraphics[width=\textwidth]{WarpingPath}
    \caption{The warping path through $D(4, 5)$ corresponding to the matching.}
  \end{subfigure}
  \caption{An alignment between two time series (left) together with its warping path (right).}
  \label{fig:warping_path}
\end{figure}

Let $s, s' \in \Series$ be two time series with lengths $m$ and $n$ respectively.
The cost $C_P(s, s')$ of a warping path in $P\in\P_{m, n}$ is given by
\[
  C_P(s, s') \define \sum_{(i, j) \in V(P)} \left( s_i - s'_j \right)^{2}.
\]
We say that a warping path \emph{aligns} elements of $s$ and $s'$,
where $s_i$ and $s'_j$ are aligned by $P$ if $(i, j) \in V(P)$.
An alignment between two time series as well as the corresponding
warping path is depicted in Figure~\ref{fig:warping_path}.

The \emph{dtw-distance} between $s$ and $s'$ is the minimum cost
of any warping path of order $m \times n$, that is,
\[
  \dtw(s, s') \define \min_{P \in \P_{m, n}} \sqrt{ C_P(s, s')}.
\]
A warping path $P$ with $\sqrt{C_P(s,s')}=\dtw(s, s')$ is called
an \emph{optimal warping path}.
An optimal warping path between two time series can
be found in quadratic time using dynamic programming~\cite{sakoe_chiba_band}. Note that the
dtw-distance is not a metric.

\paragraph*{Global Warping Path Constraints.}

The definition of warping paths between time series allows alignments
straying far from the diagonal (vertices of $D_{m,n}$ with equal
coordinates, that is, $(i,i) \in V(D_{m,n})$).  This leads to optimal
warping paths aligning elements which are far apart in time.  Since
this behavior is generally considered undesired, global constraints
restricting the set of warping paths are added in
practice~\cite{berndt1994using,chu2002iterative}.

Formally, a \emph{global constraint} is a relation
$\mathcal{R} \subseteq [m] \times [n]$ such that the preimage of $j$
under $\mathcal{R}$, defined as
$\mathcal{R}^{-1}[j] \define \{j \in [n] \mid (i, j) \in \mathcal{R} \}$
corresponds to an \emph{interval} of values, \ie to a set
\[
  [a, b] \define \{ a, \ldots ,b \} \subseteq [m]
\]
determined by endpoints $a, b \in \N$ with $a \leq b$.  A warping path
$P \in \P_{m, n}$ satisfies the global constraint $\mathcal{R}$ iff
all vertices $(i, j) \in V(P)$ are contained in $\mathcal{R}$.  The
definition of the dtw-distance is easily amended to the constrained
case.

There are two commonly used constraints regarding admissible warping paths:
The Sakoe-Chiba band~\cite{sakoe_chiba_band} of width
$r \in \N$ restricts path to vertices $v=(i, j)$ such that
$|i - j| \leq r$.
Similarly, the Itakura parallelogram~\cite{itakura_parallelogram} of slope $\sigma \ge 1$ restricts
the the warping path to vertices $(i, j)$
such that
\begin{equation}
  \frac{1}{\sigma} \leq \frac{j}{i} \leq \sigma \textrm{ and}
  \ \frac{1}{\sigma} \leq \frac{n-j+1}{m - i+1} \leq \sigma.
\end{equation}
Both restrictions have the additional advantage of decreasing the
running time needed to compute the dtw-distance.
Note that both restrictions require $s$ and $s'$ to be
compatible with respect to their lengths $m$ and $n$.  Specifically,
whenever $|m - n| > r$, no path in $\P_{m, n}$ is contained in the
Sakoe-Chiba band. The Itakura parallelogram in turn requires that
$\frac{1}{\sigma} \le \frac{m}{n} \le \sigma$ to allow for any feasible path.

\paragraph*{Fréchet Mean.}
Consider a finite sample $\mathcal{X} \define \{ s^{1}, \ldots, s^{k}\}$ of time series.
The \emph{Fréchet function} $F$ measures the
average squared dtw-distance of a time series $z \in \Series$ to the
set $\mathcal{X}$ and is defined as
\[
  F(z) \define \frac{1}{k}\sum_{i = 1}^{k} \dtw\left( z, s^{i} \right)^2.
\]
The \DTWMean problem is to find a mean time series, that is, a time
series~$z$ minimizing $F(z)$. The decision problem is defined as follows.

\problemdef{\DTWMean}
{A list of~$k$ time series $s_1,\ldots, s_k$ and~$c\in\Q$.}
{Is there a time series~$z$ such that $F(z)\le c$?}

It is known that a mean always exists (not necessarily unique)~\cite{JS16}.
In fact, there always exists a mean of length bounded linearly in the input size.

\begin{thm}[\cite{JS16}]
  \label{cor:mean_length_upper_bound}
  Let $s^{1}, \ldots, s^{k}$ be time series with lengths $m_1, \ldots, m_k$. There exists a mean $z \in \Series$
  of length at most
  \begin{equation}
    \label{eq:mean_length_upper_bound}
    \sum_{l = 1}^{k} m_{l} - 2(k - 1).
  \end{equation}
\end{thm}

Moreover, it is known that the optimal warping paths between a mean and the input time series determine the value of a mean element to be the arithmetic mean of the input values aligned to it.

\begin{lem}[\cite{SJ18}]
  \label{lem:mean_element}
  Let $z=(z_1,\ldots,z_L) \in \Series$ be a mean of the time series $s^{1}, \ldots, s^{k}$ and
  let $P_1,\ldots,P_k$ be the corresponding optimal warping paths. Then, for~$j\in[L]$, it holds that
  \begin{equation}
    z_j = \frac{\sum_{l \in [k]} \sum_{(i, j) \in V(P_l)} s^{l}_i}{\sum_{l \in [k]} |\{(i, j) \in V(P_l)\}|}.
  \end{equation}
\end{lem}

Note that the global constraints (Sakoe-Chiba band and Itakura
parallelogram) mentioned above can be added to the \DTWMean problem by
restricting the warping paths between $z$ and $s_{l}$ using global
constraints $\mathcal{R}_l \subseteq [m_l] \times [n]$ for all
$l \in [k]$.

\section{Bounding the Mean Domain}
\label{sec:bounds}

In order to obtain tight mathematical programming formulations for \DTWMean we proceed to bound not only the length of a mean but also the individual values.
In the following, we focus on upper bounds for the mean values, the case of
lower bounds is symmetric.

Let $s^{1}, \ldots, s^{k}$ be time series with lengths $m_1, \ldots, m_k$.
A simple upper bound for the value of a mean element is given by the maximum value occurring in the input
\[
  \ub^{\simple} \define \max_{l \in [k]} \max_{i \in [m_l]} s^{l}_{i}.
\]
While this bound is easily computed in linear time, it does not translate into a nontrivial bound on the Fréchet function $F$.

Regarding improvement, note that, by definition of a warping path,
every mean element is aligned with a consecutive subseries of each
input series $s^l$ defined by an interval
$I_l = [a_l,b_l] \subseteq [m_l]$. \Cref{lem:mean_element} yields the
improved bound
\begin{equation}
  \label{eq:exact_bound}
  \ub^{\improved} \define \max_{\substack{I_1, \ldots, I_k \\ I_l=[a_l,b_l]\\ 1\le a_l \le b_l \le m_l}}
  \frac{\sum_{l \in [k]} \sum_{i\in I_l}s^l_i}{\sum_{l \in [k]} |I_l|}.
\end{equation}
In order to compute~\eqref{eq:exact_bound} in polynomial time,
we follow an approach by \textcite{EH97}.
Observe that for $K \in \Q$, we have that

\begin{equation}
  \label{eq:bound_derivation}
  \begin{aligned}
    \ub^{\improved} \leq K \quad &
    \Leftrightarrow \quad
    \frac{\sum_{l \in [k]} \sum_{i\in I_l}s^l_i}{\sum_{l \in [k]} |I_l|} \leq K
    & \forall I_1, \ldots, I_k\\
    & \Leftrightarrow \quad
    \sum_{l \in [k]} \left(\sum_{i\in I_l}s^l_i - K |I_l|\right) \leq 0
    & \forall I_1, \ldots, I_k\\
    & \Leftrightarrow \quad
    \max_{\substack{I_1, \ldots, I_k \\ I_l=[a_l,b_l]\\ 1\le a_l \le b_l \le m_l}}
    \sum_{l \in [k]} \left(\sum_{i\in I_l}s^l_i - K |I_l|\right) \leq 0\\
    & \Leftrightarrow \quad
    \sum_{l \in [k]}\max_{\substack{I_l=[a_l,b_l]\\ 1\le a_l \le b_l \le m_l}}
    \left(\sum_{i\in I_l}s^l_i - K |I_l|\right) \leq 0.
  \end{aligned}
\end{equation}
Define the function $f \colon \Q \to \Q$ as follows
\begin{equation}
  \label{eq:bound_function}
  f(K) \define
  \sum_{l \in [k]}\max_{\substack{I_l=[a_l,b_l]\\ 1\le a_l \le b_l \le m_l}}
  \left(\sum_{i\in I_l}s^l_i - K |I_l|\right).
\end{equation}
Then, $\ub^{\improved} \leq K$ if and only if~$f(K)\le 0$.
Note that~$f$ is a piecewise linear decreasing function since it is a sum of~$k$ piecewise linear decreasing functions (maxima of linear decreasing functions).
Thus, computing~$\ub^{\improved}$ corresponds to finding the root of~$f$.
Note that we can evaluate~$f$ by enumerating all
intervals of all of input time series in~$\O\big(\sum_{l\in[k]}m_l^2\big)$ time.
In order to compute an (approximate) root of $f$, we
employ a binary search. Note that
\[
  \frac{1}{k} \sum_{l \in [k]} \max_{i \in m^l} s^l_i \le \ub^{\improved} \le \ub^{\simple}.
\]
We can therefore approximate $\ub^{\improved}$ in polynomial time.

Based on the upper bound $\ub^{\improved}$ and its counterpart $\lb^{\improved}$
we can bound the dtw-distance between a mean~$z$ and input series~$s^l$ as follows
\begin{equation}
  \label{eq:frechet_bound}
  \dtw(s^{l}, z)^{2} \geq
  \sum_{i \in [m_l]}
  \begin{cases}
    \left( s^{l}_i - \ub^{\improved} \right)^2 & \textrm{ if } s^{l}_i > \ub^{\improved}, \\
    \left( s^{l}_i - \lb^{\improved} \right)^2 & \textrm{ if } s^{l}_i < \lb^{\improved}, \\
    0 & \textrm{ otherwise}.
  \end{cases}
\end{equation}
This bound on $\dtw(s^{l}, z)^{2}$ translates into a nontrivial lower
bound on the Fréchet function $F$ analogously to the well-known \texttt{LB\_Keogh} bound
for the dtw-distance~\cite{keogh_bound}.

\begin{figure}[t]
  \begin{centering}
    \begin{subfigure}[t]{\textwidth}
      \includegraphics[width=\textwidth,height=.3\textwidth]{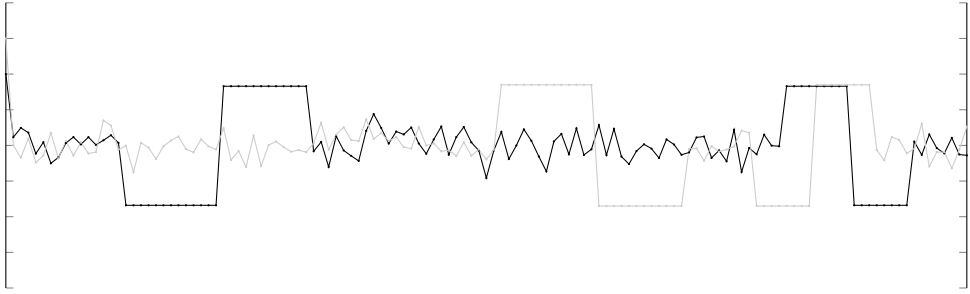}
      \caption{Two sample time series consisting of \num{129} data points each from the ``TwoPatterns'' instance from~\cite{UCR18}.}
    \end{subfigure}
    \begin{subfigure}[t]{\textwidth}
      \includegraphics[width=\textwidth,height=.3\textwidth]{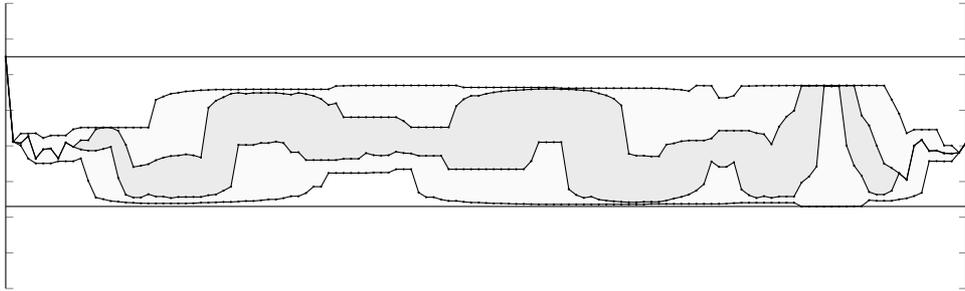}
      \caption{The improved bounds $\lb^{\improved}$ and $\ub^{\improved}$ in
        case of absence of global constraints as well as with
        respect to Itakura parallelograms of slopes \num{1.5} and
        \num{1.1} (encompassing the shaded areas).}
    \end{subfigure}
  \end{centering}
  \caption{An illustration of the improved lower and upper bounds on the mean
    values.}
  \label{fig:illustration_bound}
\end{figure}

Global constraints like the Sakoe-Chiba band or the Itakura
parallelogram (\Cref{sec:prelims}) can be combined with the upper
bound~\eqref{eq:exact_bound} by restricting the set of intervals
$I_l$.  In this case, the bound $\ub^{\improved}$ becomes dependent
upon the index~$j$ of the mean element $z_j$ under consideration.
Specifically, we let
\[
  \ub^{\improved}_j \define \max_{\substack{I_1, \ldots, I_k \\ I_l=[a_l,b_l] \subseteq \mathcal{R}_l^{-1}[j]}}
  \frac{\sum_{l \in [k]} \sum_{i\in I_l}s^l_i}{\sum_{l \in [k]} |I_l|}.
\]
Clearly, the computation of these bounds does not differ much from
the computation of the original bound $\ub^{\improved}$.
Analogously, we let $\lb^{\improved}_j$ be the corresponding lower
bound value (see \Cref{fig:illustration_bound} for an example).
The bound on the dtw-distance can be generalized to the constrained
case in a similar fashion:
\begin{equation*}
  \label{eq:constrained_frechet_bound}
  \dtw(s^{l}, z)^{2} \geq
  \sum_{i \in [m_l]}
  \min_{j \in [n]}
  \begin{cases}
    \left( s^{l}_i - \ub^{\improved}_j \right)^2 & \textrm{ if } s^{l}_i > \ub^{\improved}_j, \\
    \left( s^{l}_i - \lb^{\improved}_j \right)^2 & \textrm{ if } s^{l}_i < \lb^{\improved}_j, \\
    0 & \textrm{ otherwise}.
  \end{cases}
\end{equation*}

\section{Formulations}
\label{sec:formulations}

In the following, we will give multiple mixed integer nonlinear
programming (MINLP) formulations of the \DTWMean problem.
To this end, we consider the $k$ diagonal grid graphs
$D_1 ,\ldots, D_k$, where $D_l \define D(m_l, N)$ for $N$ being the
upper bound on the mean length given by~\eqref{eq:mean_length_upper_bound}.
Since we do not a priori know the exact length of a mean $z$, we model the mean length using
binary variables. Thus, we can solve the \DTWMean problem using a
single (albeit large) MINLP.
For notational convenience, we let
\begin{equation}
  \begin{aligned}
    V_l \define V(D_l), & \qquad  A_l \define A(D_l), \\
    s_l \define s_{m_l, N}, & \qquad T_l \define T(D_l).
  \end{aligned}
\end{equation}

\subsection{A Vertex-Based Formulation}

We begin by introducing a formulation based on binary variables
denoting whether or not a vertex in $V_l$ is part of the warping path
$P_l$ aligning $s^{l}$ and a mean~$z$. Since we do not know the length of~$z$, we include variables $x_j$ determining the length:
\begin{equation}
  \label{eq:mean_length}
  \begin{aligned}
    \sum_{j = 1}^{N} x_j = 1 \\
    x_{j} \in \{0, 1\} && \forall j \in [N]. \\
  \end{aligned}
\end{equation}
The membership of vertices from $V_l$ in $P_l$ is determined by binary
variables~$y_v^{l}$ for $v \in V_l$. It is clear that the source~$s_l$ of
$D_l$ must be contained in $P_l$ as well as one of the vertices in
$T_l$. Furthermore, if $u \in P_l$, then either $u=(m_l, j)$ for $j$ being the mean
length, or one of the out-neighbors of~$u$ must be in $P_l$ as
well.
Thus, the set of vertices of the $k$ warping paths can be described as follows:
\begin{equation}
  \label{eq:dtw_vertex_paths}
  \begin{aligned}
    y_v^{l} \in \{0, 1\}
    && \forall l \in [k], v \in V_l \\
    y_{s_l}^{l} = 1
    && \forall l \in [k] \\
    y_u^{l} \leq \smashoperator{\sum_{(u, v) \in A_l}} y_v^{l}
    && \forall l \in [k], u \in V_l \setminus T_l \\
    y_u^{l} \leq \smashoperator{\sum_{(u, v) \in A_l}} y_v^{l} + x_j
    && \forall l \in [k], u = (m_l, j) \in  T_l. \\
  \end{aligned}
\end{equation}
Lastly, the distance between input elements $s^l_i$ and mean elements $z_j$ must be included in
order to model the objective function (the Fr\'echet function~$F$). However, not all $s^l_i$ and $z_j$ are necessarily aligned in an optimal solution.
Let $d_v^{l}$ for $v = (i, j)$ be the variable denoting the cost contribution of~$v$ to the cost of~$P_l$.
We note that $d_v^{l} \geq 0$ (clearly, if~$v\not\in P_l$, then $d_v^{l} = 0$) and, more importantly, that
$d_v^{l} \leq \left( M^{l}_v \right)^2$, where
\[
  M^{l}_v \define \max (|s^{i}_l - \lb^{\improved}_j|, |\ub^{\improved}_j - s^{i}_l|).
\]
We can therefore include the mean values $z_j$ and the
corresponding distances via
\begin{equation}
  \label{eq:dtw_vertex_quad_distances}
  \begin{aligned}
    & d_v^{l} \geq (z_j - s_{i}^l)^2 - \left( M^{l}_v \right)^2 \left( 1 - y_v^{l} \right)
    & \forall l \in [k], v = (i, j) \in V_l  \\
    & 0 \leq d_v^{l} \leq \left( M^{l}_v \right)^2
    & \forall l \in [k], v \in V_l \\
    & z_j \in [\lb^{\improved}_j, \ub^{\improved}_j]
    & \forall j \in [N]. \\
  \end{aligned}
\end{equation}
The objective can then be expressed solely in terms of the
variables $d^{l}_v$, yielding the complete formulation:
\begin{equation}
  \label{eq:dtw_vertex}
  \tag{DTW-V}
  \begin{aligned}
    \min \quad & \frac{1}{k} \sum_{l \in [k]} \sum_{v \in V_l} d^{l}_v \\
    \st \quad & \eqref{eq:mean_length},
    \eqref{eq:dtw_vertex_paths}\textrm{, and }
    \eqref{eq:dtw_vertex_quad_distances}.
  \end{aligned}
\end{equation}
Overall, both the number of variables and the number of constraints
are in~$\O(k^2n^2)$, where~$n=\max (m_1,\ldots,m_k)$.  Note that the
constrains are linear except for the quadratic distance constraints
on~$d^l_v$ in~\eqref{eq:dtw_vertex_quad_distances}. It is
straightforward to replace these constraints by linear ones:
\begin{equation}
  \label{eq:dtw_vertex_linear_distances}
  \begin{aligned}
    & w_v^{l} \geq (z_j - s_{i}^l) - M^{l}_v  \left( 1 - y_v^{l} \right) & \forall l \in [k], v = (i, j) \in V_l  \\
    & w_v^{l} \geq (s_{i}^l - z_j) - M^{l}_v  \left( 1 - y_v^{l} \right) & \forall l \in [k], v = (i, j) \in V_l  \\
    & 0 \leq w_v^{l} \leq M^{l}_v & \forall l \in [k], v \in V_l \\
    & z_j \in [\lb^{\improved}_j, \ub^{\improved}_j] & \forall j \in [N].
  \end{aligned}
\end{equation}
and minimize
$\sum_{l \in [k]} \sum_{v \in V_l} \left( w_v^{l} \right)^2$ instead, thus
modeling \DTWMean as a mixed-integer quadratic program (MIQP).  In any
case, the distance constraints are of the so-called \emph{big-$M$} type,
known to be numerically more challenging and to yield poor
lower bounds in general.

\begin{remark}[Size]
  \label{rem:size}
  A notable disadvantage of
  formulation~\eqref{eq:dtw_vertex_linear_distances} is its size in
  terms of number of variables: For a set of $k$ time series of a uniform length
  of $m$, formulation~\eqref{eq:dtw_vertex_linear_distances} consists of
  $\Theta(km^2)$ many variables, making it challenging to solve.
  The large size stems from the fact that both the mean
  length and the alignments are entirely unknown beforehand.

  Interestingly, the addition of global constraints
  alleviates both problems, greatly facilitating the practical
  tractability of the problem: Firstly, the Sakoe-Chiba band and the
  Itakura parallelogram restrict the maximum mean length~$N$ to be relatively
  close to the input length~$m$. Secondly, many alignments are excluded beforehand,
  eliminating the corresponding variables and constraints.
  Thus, the addition of global constraints is not only advantageous
  with respect to qualitative considerations, but also makes the
  problem computationally more tractable.
\end{remark}

\subsection{An Arc-Based Formulation}
Recall that the $y$-variables in the vertex-based
formulation~\eqref{eq:dtw_vertex} model warping paths through the
digraphs $D_l$. Conventionally, paths through networks are described
in terms of unit network flows (see~\cite[pp. 173]{bible}). We will
therefore proceed to give an arc-based formulation in addition to the
vertex-based formulation introduced above. We will then discuss the merits of this formulation
(in \Cref{sec:computational} we also conduct several computational experiments).

Formally, we formulate the problem of finding a set of $k$ warping
paths through the graphs $D_l$ in terms of binary arc variables
forming a set of unit flows, \ie
\begin{equation}
  \label{eq:flow_conservation}
  \begin{aligned}
    &f^{l}_{a} \in \{0, 1\} &\forall l \in [k], a \in A_l, \\
    & f^{l}(\delta^{+}(v)) - f^{l}(\delta^{-}(v)) =
    \begin{cases}
      1 & \textrm{ if } v = s_l, \\
      -x_j & \textrm{ if } v = (m_l, j) \in T_l,\\
      0 & \textrm{ otherwise.}
    \end{cases} &\forall l \in [k], v \in V_l.
  \end{aligned}
\end{equation}
where $\delta^{+}(u)$ and $\delta^{-}(u)$ denote the outgoing and
incoming arcs of a vertex $u$ respectively.  This definition of
warping paths in terms of flows enables us to express the
variables~$y_{v}^{l}$ in terms of the corresponding flow variables
by means of the following coupling constraints
\begin{equation}
  \label{eq:flow_coupling}
  \begin{aligned}
  &y_v^{l} =
  \begin{cases}
    f^{l}(\delta^{+}(v)) + x_j & \textrm{if } v = (m_l, j) \in  T_l, \\
    f^{l}(\delta^{+}(v))  & \textrm{otherwise}
  \end{cases} & \forall l\in [k], v\in V_l.
  \end{aligned}
\end{equation}
The resulting flow-based formulation of \DTWMean is given as
\begin{equation}
  \label{eq:dtw_arc}
  \tag{DTW-A}
  \begin{aligned}
    \min \quad & \frac{1}{k} \sum_{l \in [k]} \sum_{v \in V_l} d^{l}_v \\
    \st \quad &
    \eqref{eq:mean_length}, \eqref{eq:dtw_vertex_quad_distances},
    \eqref{eq:flow_conservation}, \textrm{ and }
    \eqref{eq:flow_coupling}.
  \end{aligned}
\end{equation}
In practice, we can use the constraints~\eqref{eq:flow_coupling} in order
to entirely eliminate the $y$-variables.

As for the differences between the formulations: While the digraphs $D_l$
are relatively sparse, we still roughly triple the number of variables
required to model all warping paths in \eqref{eq:dtw_arc}. Since the size of the formulations
is already significant (see Remark~\ref{rem:size}), a further increase
in size seems undesirable.

However, we can also judge different formulations in terms of their
tightness. Specifically, let $\NLP(\cdot)$ be the \emph{NLP
  relaxation} of a formulation, \ie the nonlinear program (NLP)
obtained by dropping integrality requirements from a MINLP.  A
formulation $A$ is said to be \emph{tighter} than $B$ iff
$\NLP(A) \subseteq \NLP(B)$. A tighter formulation yields stronger
bounds, making branch-and-bound procedures more efficient.  It is easy to
show that \eqref{eq:dtw_arc} is tighter than \eqref{eq:dtw_vertex}:
\begin{lem}
  $\proj_{y, d} \left( \NLP \eqref{eq:dtw_arc} \right) \subseteq \NLP \eqref{eq:dtw_vertex}.$
\end{lem}
\begin{proof}
  Let $(x, f, y ,d)$ be a solution of $\NLP \eqref{eq:dtw_arc}$, \ie
  satisfying all constraints of~\eqref{eq:dtw_arc} except for the
  integrality condition $f^{l}_{a} \in \{0, 1\}$ which is relaxed to
  $0 \leq f^{l}_a \leq 1$ for all $l \in [k]$, $a \in A_l$. In order
  to prove $(x, y, d) \eqref{eq:dtw_vertex}$ it is sufficient to show
  that $y$ satisfies the relaxation of the
  constraints~\eqref{eq:dtw_vertex_paths}.

  Consider a vertex $u \in V_l$ for some $l \in [k]$. Based on the
  coupling constraints~\eqref{eq:flow_coupling}, $y_{u}^l$ must be non-negative.
  The flow constraints~\eqref{eq:flow_conservation} in turn imply that $y_{u}^{l}$
  is less than or equal to one, where equality holds for $y_{s_l}^l$.
  If $u \notin T_l$ it holds that
  \begin{equation}
    \begin{aligned}
      y_u^{l}
      = & f^{l}(\delta^{+}(u)) = \sum_{(u, v) \in A_{l}} f^{l}(u, v)
      \leq \sum_{(u, v) \in A_{l}} f^{l}(\delta^{-}(v)) \\
      = & \sum_{(u, v) \in A_l : v \in V_l \setminus T_l} f^{l}(\delta^{-}(v))
      + \sum_{(u, v) \in A_l : v = (m_l, j) \in T_l} f^{l}(\delta^{-}(v)) + x_j \\
      = & \sum_{(u, v) \in A_{l}} y_v.
    \end{aligned}
  \end{equation}
  The case of $u \in T_l$ can be treated analogously.
\end{proof}

\subsection{Distance Formulations}

Recall that the distance constraints involving the variables $d_{v}^{l}$
are given by
\begin{equation}
  \label{eq:distance_bigM}
  \begin{aligned}
    d_v^{l} \geq (z_j - s_{i}^l)^2 - \left( M^{l}_v \right)^2 \left( 1 - y_v^{l} \right)
    && \forall l \in [k], v=(i,j) \in V_l,
  \end{aligned}
\end{equation}
where $y_{v}^{l}$ denotes whether or not a vertex $v$ is contained in
a warping path $P_l$. Determining the optimal warping paths, while
constituting the key difference between the arc-based~\eqref{eq:dtw_arc} and vertex-based~\eqref{eq:dtw_vertex}
formulations above, is independent of how the distances are modeled exactly.
We can therefore study different formulations of the objective independently of the underlying graph
model.
As mentioned before, the constraints~\eqref{eq:distance_bigM} are
big~$M$ constraints switched on and off by the $y$-variables.
Since these constraints generally yield poor relaxations, it is worth
investigating alternative modeling techniques.

\paragraph{A Perspective Reformulation.}

In the following, we will derive an alternative formulation of the
distance constraints in order to avoid the big-$M$ constraints present
in the original formulation. To this end, we will consider a fixed
index $l \in [k]$ and vertex $v \in V_l$.  A straightforward
reformulation of the distance constraints can be obtained by weighing
the quadratic distance between mean elements and input elements with
the $y$-variables, \ie requiring that
\begin{equation*}
  \label{eq:distance_cubic}
  \begin{aligned}
    d_v^{l} \geq  y_v^{l} \cdot (z_j - s_{i}^l)^2
    && \forall l \in [k], v=(i,j) \in V_l. \\
  \end{aligned}
\end{equation*}
Unfortunately, these inequalities are non-convex, making it
practically impossible to solve the fractional relaxations.
We want the variables $d_{v}^{l}$, $z_j$, and $y_{v}^{l}$ variables to
be contained in the union of the following convex bounded sets
\begin{equation*}
  \begin{aligned}
    P^{0} &\define \{ (d_v^ {l}, z_j, y_v^{l}) \mid y_v^{l} = 0, d_v^{l} = 0, z_j \in [\lb^{\improved}_j, \ub^{\improved}_j] \}, \textrm{ and} \\
    P^{1} &\define \{ (d_v^{l}, z_j, y_v^{l})  \mid y_v^{l} = 1, \left( M^{l}_v \right)^2 \geq d_v^{l} \geq (z_j - s_{i}^l)^2, z_j \in [\lb^{\improved}_j, \ub^{\improved}_j] \}.
  \end{aligned}
\end{equation*}
In order to obtain a convex optimization problem we would like to have a
description of $\conv(P^{0} \cup P^{1})$ in terms of a set of convex inequalities.
To this end, we can use the so-called \emph{perspective reformulation}~\cite{perspective_reformulation}.
The \emph{perspective function} of a function $f : \R^{n} \to \R^{m}$ is the function
$\tilde{f} : \Rp \times \R^{n} \to \R^{m}$ defined as
\[
  \tilde{f}(\lambda, x) \define
  \begin{cases}
    \lambda \cdot f (x / \lambda) & \textrm{ if } \lambda > 0, \textrm{ and} \\
    0 & \textrm{ if }  \lambda = 0. \\
  \end{cases}
\]
It holds that the perspective of a convex function is also convex.
The key observation regarding the perspective function is the following:
\begin{thm}[\cite{perspective_reformulation, convex_disjunctive}]
  \label{thm:persp}
  Let $f^{t} : \R^{n} \to \R^{m_t}$ for $t \in T$ be a set of functions such that
  the sets
  \[
    K^{t} \define \{x \in \R^{n} \mid f^{t}(x) \leq 0 \}
  \]
  are convex and bounded. Then, $x \in \conv ( \bigcup_{t \in T} K^{t})$ if and only if
  \[
      x = \sum_{t \in T} x^{t}, \quad \sum_{t \in T} \lambda_t = 1, \quad \tilde{f}^{t} (\lambda_{t}, x^{t}) \leq 0, \quad \lambda_t \geq 0 \: \forall t \in T.
  \]
\end{thm}
Based on \Cref{thm:persp}, we obtain the desired
description including an additional variable $\overline{z}_j$:
\begin{equation*}
  \begin{aligned}
    z_j- \overline{z}_j & \in \left[ (1 - y_{v}^{l}) \cdot \lb^{\improved}_j, (1 - y_{v}^{l}) \cdot \ub^{\improved}_j \right] \\
    d_{v}^{l} & \leq y_{v}^{l} \cdot \left( M^{l}_v \right)^2 \\
    \overline{z}_j & \in \left[ y_{v}^{l} \cdot \lb^{\improved}_j, y_{v}^{l} \cdot \ub^{\improved}_j \right] \\
    d_{v}^{l} & \geq \left( s_{i}^{l}  - \overline{z}_j/y_{v}^{l}  \right)^2, \textrm{ if } y_{v}^{l} > 0.
  \end{aligned}
\end{equation*}
Note that the domain of the last inequality cannot easily be extended
to include $y_{v}^{l} = 0$, which is due to the piecewise definition
of the perspective function. Observe that for $y_{v}^{l} = 0$ the
other inequalities already imply that
$d_{v}^{l} = \overline{z}_j = 0$.  Still, it is
well-known~\cite{perspective_reformulation} that solvers frequently
struggle with numerical problems when encountering perspective
functions.  We therefore propose to begin instead by applying an
\emph{outer approximation}~\cite{outer_approximation} to the quadratic
constraint $d_v^{l} \geq (z_j - s_{i}^l)^2$.  The outer approximation
of a set $S_f \define \{x \in \R^{n} \mid f(x) \leq 0\}$ given in terms of a convex
differentiable function $f : \R^{n} \to \R$ is given as the polyhedron
defined by
\[
  \{x \in \R^{n} \mid f(x^{i}) + \langle \nabla f(x^{i}), x - x^{i} \rangle \leq 0 \: \forall i \in [k]\}
  \supseteq S_f
\]
based on a set $\{x^{1}, \ldots, x^{k}\} \subseteq \R^{n}$ of supporting points.
\begin{figure}
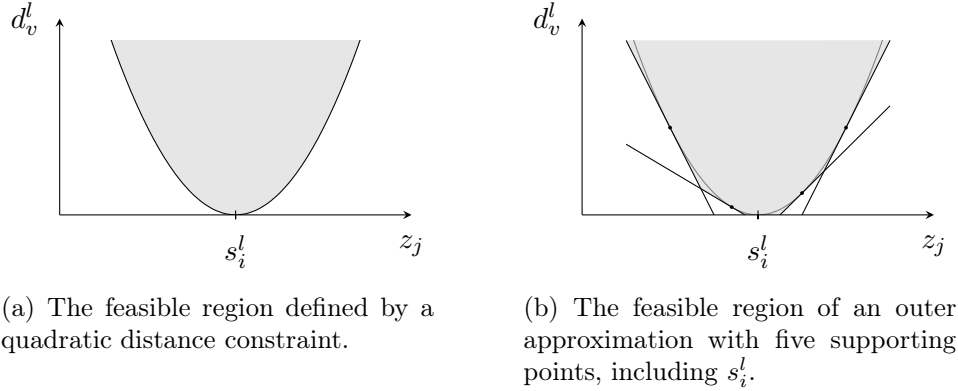

  \begin{centering}
    \begin{subfigure}[t]{0.45\textwidth}
      \includegraphics[width=\textwidth]{OuterOriginal}
      \caption{The feasible region defined by a quadratic distance constraint.}
    \end{subfigure}
    \hspace{1cm}
    \begin{subfigure}[t]{0.45\textwidth}
      \includegraphics[width=\textwidth]{OuterApprox}
      \caption{The feasible region of an outer approximation with five supporting points,
        including $s_{i}^{l}$.}
    \end{subfigure}
  \end{centering}
  \caption{An illustration of the outer approximation of a quadratic distance constraint.}
  \label{fig:outer_approx}
\end{figure}
In our case, a fixed value of $\overline{z}_j$ yields the inequality
\[
  (\overline{z}_j)^2 - (s_{i}^l)^2 \geq 2 (\overline{z}_j - s_{i}^l) z_j - d_v^{l}.
\]
See \Cref{fig:outer_approx} for an example of such an outer
approximation. Note that this inequality does not depend on the value
$\overline{d}_{v}^{l}$ corresponding to~$\overline{z}_j$.  For a set
$\{z^1_j, \ldots, z^{r}_j\}\subseteq \R$ of supporting points, we
obtain a linear system of inequalities $A_r \cdot (z_j, d_{v}^{l}) \leq b_r$, where
\begin{equation*}
  \begin{aligned}
    A_r \define
    \begin{pmatrix}
      2(z^1_j - s_{i}^{l}) & -1 \\
      \vdots & \vdots \\
      2(z^r_j - s_{i}^{l}) & -1 \\
    \end{pmatrix}, \quad &
    b_r \define
    \begin{pmatrix}
      (z^{1}_{j})^2 - (s_{i}^l)^2 \\
      \vdots \\
      (z^{r}_{j})^2 - (s_{i}^l)^2
    \end{pmatrix}.
  \end{aligned}
\end{equation*}
We proceed to apply the perspective reformulation to the sets $P^{0}$ and
\begin{equation*}
\begin{aligned}
  Q^{1} \define \{ (d_v^{l}, z_j) \mid \: & y_{v}^{l} = 1,
  d_v^{l} \leq \left( M^{l}_v \right)^2, A_r \cdot (z_j, d_{v}^{l}) \leq b_r, \\
  & z_j \in [\lb^{\improved}_j, \ub^{\improved}_j] \}.
\end{aligned}
\end{equation*}
Note that since both $P^{0}$ and $Q^{1}$ are polytopes, the set
$\conv(P^{0} \cup Q^{1})$ is a polytope as well, alleviating the
numerical problems of the perspective reformulation based on nonlinear
functions. Indeed, optimizing over the union of polytopes is known as
a \emph{disjunctive programming problem} (see~\cite{disjunctive}). In our case we
obtain the following constraints:
\begin{equation*}
  \begin{aligned}
    z_j- \overline{z}_j &\in \left[ (1 - y_{v}^{l}) \cdot \lb^{\improved}_j, (1 - y_{v}^{l}) \cdot \ub^{\improved}_j \right] \\
    d_{v}^{l} &\leq y_{v}^{l} \cdot \left( M^{l}_v \right)^2 \\
    \overline{z}_j &\in \left[ y_{v}^{l} \cdot \lb^{\improved}_j, y_{v}^{l} \cdot \ub^{\improved}_j \right] \\
    A_r \cdot (\overline{z}_j, d_{v}^{l}) &\leq y_{v}^{l} \cdot b_r \\
  \end{aligned}
\end{equation*}
In order to solve the \DTWMean problem we still have to settle on a
set of supporting points: If the set is too small or unevenly spaced,
then the outer approximation is not sufficiently tight, and the
error with respect to the actual quadratic function becomes too large.
On the other hand, each support point increases the size of
the resulting program, slowing down the solution process.

To avoid the problem of having to select a suitable set beforehand, we
\emph{separate} inequalities as needed: We begin with a support set
consisting only of $s_{i}^{l}$ and solve the resulting problem,
obtaining a solution consisting in part of values for $d_v^{l}$ and
$\overline{z}_j$. If these values violate the quadratic constraint
sufficiently much, we add the value of $\overline{z}_j$ to the
supporting points (thereby cutting off the solution) and resolve.
In practice, all state-of-the-art MINLP solvers offer so-called
\emph{callback functions} in order to support the separation of
additional constraints during the solution process.

Still, the approach comes at a price in terms of size: For each
distance variable~$d_v^l$, we need one additional variable and several
additional constraints.

\paragraph{Implicit Distances.}

In the following we will consider the framework introduced
by~\textcite*{unified_framework}. The authors observed that the
big-$M$ modeling approach can be seen as a regularization of logical
constraints in a two-stage problem consisting of an outer binary and
an inner continuous part.
The problem can be dualized, yielding a
convex inner function which can be used to derive valid inequalities
based on outer approximations. We will adapt this approach to our
formulations.

Consider a feasible solution $(x, y)$ of
the system comprised of both~\eqref{eq:mean_length}
and~\eqref{eq:dtw_vertex_paths}, that is, a set of warping paths in
the digraphs $D_l$ connecting their respective sources with vertices
corresponding to a fixed mean length.  Let $P_l$ be the path such that
$y_{v}^{l} = 1$ if and only if $v \in V(P_l)$.  For convenience, we
let $\overline{y}_{v}^{l} \define 1 - y_{v}^{l}$ be the inverse of the
variable~$y_{v}^{l}$ denoting the absence of a vertex $v$ from
$V(P_l)$.  We rewrite the distance
constraints~\eqref{eq:dtw_vertex_quad_distances} by introducing
additional variables $\Delta_{v}^{l}$:
\begin{equation}
  \label{eq:inner_constraints}
  \begin{aligned}
    d_{v}^{l} \geq (z_j - s_{i}^l)^2 - \Delta_{v}^{l} && \forall l \in [k], v = (i, j) \in V_l  \\
    d_{v}^{l}, \Delta_{v}^{l} \geq 0 && \forall l \in [k], v = (i, j) \in V_l  \\
  \end{aligned}
\end{equation}
together with the following logical constraint:
\begin{equation}
  \label{eq:inner_logical}
  \begin{aligned}
    \Delta_{v}^{l} = 0  \enspace \textrm{if} \enspace \overline{y}_{v}^{l} = 0 && \forall l \in [k], v = (i, j) \in V_l.  \\
  \end{aligned}
\end{equation}
The problem of finding a mean $z$ corresponding to the solution~$(x,y)$ is then given as
\begin{equation}
  \label{eq:inner_problem}
  f(y) \define
  \begin{cases}
    \begin{aligned}
      \min_{d, \Delta, z} \quad & \frac{1}{k} \sum_{l \in [k]} \sum_{v \in V_l} d^{l}_v \\
      \st \quad &
      (d, \Delta, z) \textrm{ satisfy }
      \eqref{eq:inner_constraints} \textrm{ and } \eqref{eq:inner_logical},
    \end{aligned}
  \end{cases}
\end{equation}
where the last constraint ensures that the binary values
$\overline{y}$ control the range of the variables $\Delta$. If a
vertex $v$ is contained in $P_l$, \ie $\overline{y}_{v}^{l} = 0$, then
$\Delta_{v}^{l}$ is fixed to zero, which may force $d_{v}^{l}$ to a
positive value in order to satisfy the first constraint. Conversely,
if $\overline{y}_{v}^{l} = 1$, then we can set $d_{v}^{l}$ to zero save
costs. Thus, an optimal solution $(d^{*}, \Delta^{*}, z^{*})$
of~\eqref{eq:inner_problem} is given by
\begin{equation}
  \begin{aligned}
    z^{*}_j &\define \frac{\sum_{l \in [k]} \sum_{v = (i, j) \in V_l} s_{i}^{l} \cdot y_{v}^{l} }{\sum_{l \in [k]} \sum_{v = (i, j) \in V_l} y_{v}^{l}},\\
    \left( d^{*} \right)_{v}^{l} &\define
    \begin{cases}
      (z^{*}_j) - s_{i}^{l} & \textrm{ if } \overline{y}_{v}^{l} = 0, \textrm{ and} \\
      0 & \textrm{ otherwise,}
    \end{cases} \\
    \left( \Delta^{*} \right)_{v}^{l} &\define
    \begin{cases}
      0 & \textrm{ if } \overline{y}_{v}^{l} = 0, \textrm{ and} \\
      (z^{*}_j) - s_{i}^{l}  & \textrm{ otherwise.}
    \end{cases}
  \end{aligned}
\end{equation}
In order to use the techniques from~\cite{unified_framework},
we add a regularization term
$ \Omega(\Delta) \define \sum_{l \in [k]} \sum_{v \in V_l} \Omega_{v}^{l}(\Delta_{v}^{l})$,
where
\begin{equation*}
  \Omega_{v}^{l}(\alpha) \define
  \begin{cases}
    0 & \textrm{ if } |\alpha| \leq (M_{v}^{l})^2, \textrm{ and} \\
    \infty & \textrm{otherwise.}
  \end{cases}
\end{equation*}
Furthermore, we introduce a function $g(d, \Delta, z)$ encompassing the objective
and parts of the constraints of~\eqref{eq:inner_problem}:
\begin{equation*}
  g(d, \Delta, z) \define
  \begin{cases}
    \frac{1}{k} \sum_{l \in [k]} \sum_{v \in V_l} d^{l}_v,
    & \textrm{ if } (d, \Delta, z) \textrm{ satisfy } \eqref{eq:inner_constraints}, \textrm{ and}  \\
    \infty
    & \textrm{ otherwise.}
  \end{cases}
\end{equation*}
Thus, the problem becomes
\begin{equation*}
  \begin{aligned}
    \min_{d, \Delta, z} \quad & g(d, \delta, z) + \Omega(\Delta) \\
    \st \quad & \Delta \textrm{ satisfies } \eqref{eq:inner_logical},
  \end{aligned}
\end{equation*}
which can, according to~\cite[Theorem~1]{unified_framework}, be transformed
into the following saddle-point problem involving additional variables
$\alpha_v^{l}$:
\begin{equation*}
  \label{eq:saddle_point}
  \begin{aligned}[t]
    \max_{\alpha}
  \end{aligned} \quad
  \begin{aligned}[t]
    \min_{d, \delta, z} \quad
    & \sum_{l \in [k]} \sum_{v \in V_l} \frac{1}{k} d_v^{l}
    - \alpha_{v}^{l} \delta_{v}^{l}
    - \overline{y}_v^{l} |\alpha_{v}^{l}|  \left( M_{v}^{l}  \right)^2 \\
    \st \quad & (d, \delta, z) \textrm{ satisfy } \eqref{eq:inner_constraints}.
  \end{aligned}
\end{equation*}
Note that the variables $\Delta$ disappear as $\Omega(\cdot)$ is replaced by
its Fenchel conjugate.  In order to solve the saddle point problem, we can make
several observations regarding the choice of variables $\alpha$: If
$\alpha_v^{l} > 0$, then the inner problem becomes unbounded since the
objective value is strictly decreasing along increasing values of
$\delta_v^{l}$.  Thus, we can assume that $\alpha_v^{l} \leq 0$. We
can therefore define $\beta_v^{l} \define - \alpha_v^{l}$ and rewrite
the problem as:
\begin{equation}
  \label{eq:saddle_point_nonneg}
  \begin{aligned}[t]
    \max_{\beta \geq 0}
  \end{aligned} \quad
  \begin{aligned}[t]
    \min_{d, \delta, z} \quad
    & \sum_{l \in [k]} \sum_{v \in V_l}
    \frac{1}{k} d_v^{l} + \beta_{v}^{l}
    \left( \delta_{v}^{l} - \overline{y}_v^{l} \left( M_{v}^{l}  \right)^2 \right) \\
    \st \quad & (d, \delta, z) \textrm{ satisfy } \eqref{eq:inner_constraints}.
  \end{aligned}
\end{equation}
Since we already know that the objective values of~\eqref{eq:inner_problem} and
\eqref{eq:saddle_point_nonneg} must coincide, it only remains to find
suitable values of $\beta$. Specifically, if we let
\begin{equation*}
  (\beta^{*})_{v}^{l} \define
  \begin{cases}
    \frac{1}{k} & \textrm{ if } \overline{y}_{v}^{l} = 0, \textrm{ and} \\
    0 & \textrm{ otherwise},
  \end{cases}
\end{equation*}
then $(d^{*}, \delta^{*} = \Delta^{*}, z^{*})$ is a solution of the
inner optimization problem of~\eqref{eq:saddle_point_nonneg} having the
same objective value as the original~\eqref{eq:inner_problem}.

For this value of $\beta$, we can derive a cutting
plane based on the subgradients of the convex function $f$.
One subgradient of $f$ is given by
\begin{equation*}
  \left( \nabla f(y) \right)^{l}_{v} \define
  \begin{cases}
    \frac{(M_{v}^{l})^2}{k} & \textrm{ if } \overline{y}_{v}^{l} = 0, \textrm{ and} \\
    0 & \textrm{ otherwise}.
  \end{cases}
\end{equation*}
Based on this subgradient, from each feasible solution $\hat{y}$ of~\eqref{eq:dtw_vertex_paths},
we obtain a linear inequality
$f(y) \geq f(\hat{y}) + \langle \nabla f(\hat{y}), (y - \hat{y})
\rangle$ as
\begin{equation*}
  f(y) \geq f(\hat{y}) + \sum_{l \in [k]} \sum_{v \in V_l : \hat{y}_{v}^{l} = 1} \frac{(M_{v}^{l})^2}{k} \left( {y}_{v}^{l} - 1 \right),
\end{equation*}
where $f(y)$ corresponds to the value of the mean derived from the solution $y$.
To embed this approach into our formulations, we introduce an additional variable
$\eta$ denoting the objective value and require
\begin{equation}
  \label{eq:implicit_distance_constraints}
  \begin{aligned}
    \eta \geq f(\hat{y}) + \sum_{l \in [k]} \sum_{v \in V_l : \hat{y}_{v}^{l} = 1} \frac{(M_{v}^{l})^2}{k} \left( {y}_{v}^{l} - 1 \right)
    && \forall \: \hat{y} \: \textrm{satisfying } \eqref{eq:dtw_vertex_paths}
  \end{aligned}.
\end{equation}
Thus, we can reformulate~\eqref{eq:dtw_vertex} as
\begin{equation*}
  \begin{aligned}
    \min \quad & \eta \\
    \st \quad & \eqref{eq:mean_length},
    \eqref{eq:dtw_vertex_paths} \textrm{, and }
    \eqref{eq:implicit_distance_constraints}.
  \end{aligned}
\end{equation*}
The arc-based formulation~\eqref{eq:dtw_arc} can be adapted in much the same way.
We would like to point out that this reformulation of the distance constraints is
much smaller, since no variables apart from $x$, $y$, and $z$ are required.
The inequalities~\eqref{eq:implicit_distance_constraints} can be separated
whenever a feasible solution is obtained throughout the search in a
branch-and-bound tree.

\section{Computational Results}
\label{sec:computational}

All experiments were conducted using an implementation in the
\texttt{C++} programming language compiled using the \textsc{GNU}
\texttt{C++} compiler with the optimizing option \texttt{-O2}. We used
version 6.0.2 of the \textsc{SCIP}~\cite{SCIP} optimization suite and
version 8.1 of \textsc{Gurobi}~\cite{Gurobi} as underlying LP solver.
All measurements were taken on an Intel Core i7-965 processor clocked
at \SI{3.2}{\giga \hertz}.

We begin by comparing the formulations across several small instances,
generated from the ``FiftyWords'' data set of the UCR
archive~\cite{UCR18}. Specifically, we sampled sets of
$k \in \{2, 5\}$ time series (of original length~$N = 270$). Each
sampled time series was reduced to a uniform length of
$m \in \{10, 20\}$ by averaging successive disjoint blocks of
$\floor{N/m}$ values.

We measure the quality of the formulations based on the remaining gap
after a time limit of one hour has expired. The gap is given as
$(p - d)/d$, where~$p$ is the value of the best-known feasible
solution, and~$d$ is the \emph{dual bound} obtained as the minimal
relaxation value across the leaves of the branch-and-bound tree.
Measuring the gap provides a good overview over the practical
performance of the different formulations, since both the solution
times of the relaxations and the dual bounds provided by them
influence the resulting gap.  To reduce the effect of the random
sampling used to generate the instances, we measured the average
remaining gap over ten instances for each variant.  Furthermore, we
included global constraints given by both a \emph{wide}
($\sigma = \num{1.5}$) and a \emph{narrow} ($\sigma=1.1$) Itakura
parallelogram.
\begin{table}[h]
  \centering
  \caption{Remaining gap after one hour of computation for different
    instances and formulations. The distance formulations are
    denoted as \textbf{Q}uadratic, \textbf{P}erspective and
    \textbf{I}mplicit respectively.}
  \label{table:formulations}
  \begin{tabular}{cS[table-format=2.2]S[table-format=2.2]S[table-format=1.]S[table-format=3.2]S[table-format=4.2]S[table-format=1]}
  \toprule
  \multirow{2}{*}{\textbf{Variant}} & \multicolumn{3}{c}{\textbf{Arc-based}} & \multicolumn{3}{c}{\textbf{Vertex-based}} \\
  \cmidrule(l{.5cm}r{.5cm}){2-4}
  \cmidrule(l{.5cm}r{.5cm}){5-7}
  & {\textbf{Q}} & {\textbf{P}} & {\textbf{I}} & {\textbf{Q}} & {\textbf{P}} & {\textbf{I}} \\
  \midrule
  $m = 10$, $k = 2$ \\
  {free}   & 2.10 & 2.56 & {$\infty$} & 3.05 & 4.89 & {$\infty$} \\
  {wide}   & 0    & 0     & 0         & 0    & 0 & 0 \\
  {narrow} & 0    & 0     & 0         & 0    & 0 & 0 \\
  \midrule
  $m = 20$, $k = 2$ \\
  {free}   & 21.79 & 69.46 & {$\infty$} & 369.29 & 3001.95    & {$\infty$} \\
  {wide}   & 0.02  & 0.03  & {$\infty$} & 0.36   & 0.38 & {$\infty$} \\
  {narrow} & 0     & 0     & 0          & 0      & 0    & 0 \\
  \midrule
  $m = 10$, $k = 5$ \\
  {free}   & 3.75 & 15.84 & {$\infty$} & 7.92 & {$\infty$}   & {$\infty$} \\
  {wide}   & 0    & 0 & {$\infty$} & 0 & 0    & {$\infty$} \\
  {narrow} & 0    & 0 & 0 & 0 & 0    & 0 \\
  \midrule
  $m = 20$, $k = 5$ \\
  {free}   & 15.48 & {$\infty$} & {$\infty$} & 25.93 & {$\infty$}    & {$\infty$} \\
  {wide}   & 0.44  & 1.78 & {$\infty$} & 0.88 & 68.40    & {$\infty$} \\
  {narrow} & 0     & 0 & 0 & 0 & 0    & 0 \\
  \bottomrule
\end{tabular}
\end{table}
From the results, displayed in \Cref{table:formulations}, we can make
several observations: Firstly, the implicit distances clearly perform
worst. This seems to be largely due to their poor relaxation
values. Indeed, the
inequalities~\eqref{eq:implicit_distance_constraints} separated while
traversing the branch-and-bound tree are often insufficient to obtain
nontrivial lower bounds. The quadratic distances performed best
for all instances. Secondly, the arc-based
formulation~\eqref{eq:dtw_arc} performs considerably better than its
vertex-based counterpart~\eqref{eq:dtw_vertex}. Apparently, the
tighter approximation of the arc-based formulation more than
compensates for the increase in size.

Unsurprisingly, an increased size with respect to both the number $k$
of time series and the length $m$ of the time series results in larger gaps remaining after the time limit. The addition of global constraints greatly improved the
solution process, presumably because of the reduction in problem size
in terms of the number of variables and constraints.

\section{Conclusion and Future Works}
\label{sec:conclusion}

In this paper, we gave the first formulation of the \DTWMean problem as a nonlinear
optimization problem. We derived nontrivial bounds
on the mean domain, translating into lower bounds on the value of the
Fréchet function also taking into account global constraints such as the
Skaoe-Chiba band or the Itakura parallelogram.

We introduced several different nonlinear programming formulations of
\DTWMean, based on different modeling approaches to the
combinatorial structure as well as the nonlinear cost function.

We compared these formulations with respect to their computational
efficiency, measured in terms of the remaining gap after one hour
of computation, concluding that a quadratic big $M$ distance formulation
together with an arc-based model for the warping paths performs
best in practice.

Unfortunately, solving the \DTWMean problem on large-scale instances still seems out of reach.
This is likely due to the fact that the introduced
MINLP formulations, while being significant in size, yield poor lower
bounds, resulting in enormous gaps relative to the primal solutions
obtained throughout the course of optimization.  As a result, few if
any branches of the branch-and-bound tree are can be discarded, and
most of the feasible solutions have to be enumerated.  On the other hand, it
is straightforward to include global constraints into the different
MINLP formulations. The formulations can then take advantage of the
reduction in combinatorial complexity and solve the resulting problems
more efficiently.

There are several directions in which this work can be extended.  On
the one hand, the formulations introduced here yield relaxations of
insufficient quality. In order to strengthen the formulations, it might
be necessary to derive families of valid inequalities.  On the other
hand, any a priori bounds on the length of the mean series would aid
computations. Conversely, the inclusion of relaxation-based heuristics
tailored specifically to the \DTWMean problem could increase the
practical performance of the formulations as well.

\printbibliography

\end{document}